\documentclass[12pt]{article}
\usepackage{mathrsfs}
\usepackage{bm, amscd,  amsfonts, amssymb, cite,latexsym}
\usepackage{amsmath, epsfig, cite, color, colordvi, longtable}

\makeatletter \@addtoreset{equation}{section} \makeatother

\newtheorem{thm}{Theorem}[section]

\newtheorem{lem}[thm]{Lemma}
\newtheorem{exam}[thm]{Example}

\def\mid{{\,|\,}}
\def\pf{\noindent {\it Proof.\ }}
\def\qed{\hfill \rule{4pt}{7pt}}
\begin{document}
\allowdisplaybreaks

\begin{center}
{\large \bf   Congruences of Multipartition Functions Modulo Powers of Primes}

\vskip 3mm

William Y.C. Chen$^1$, Daniel K. Du$^2$, Qing-Hu Hou$^3$ and Lisa H.Sun$^4$\\[5pt]
Center for Combinatorics, LPMC-TJKLC \\
Nankai University, Tianjin 300071, P. R. China

\vskip 3mm

 E-mail: $^1$chen@nankai.edu.cn,$^2$dukang@mail.nankai.edu.cn, $^3$hou@nankai.edu.cn,
 $^4$sunhui@nankai.edu.cn
\end{center}

\noindent
{\bf Abstract.}
Let $p_r(n)$ denote the number of $r$-component multipartitions of $n$, and let $S_{\gamma,\lambda}$ be the space spanned by  $\eta(24z)^\gamma \phi(24z)$, where $\eta(z)$ is the Dedekind's eta function and $\phi(z)$ is a holomorphic modular form in $M_\lambda({\rm SL}_2(\mathbb{Z}))$.
In this paper, we show that the generating function of $p_r(\frac{m^k n +r}{24})$ with respect to $n$ is congruent to a function in the space $S_{\gamma,\lambda}$ modulo $m^k$. As special cases, this   relation  leads to many well known congruences including the Ramanujan congruences of $p(n)$ modulo $5,7,11$ and Gandhi's congruences of $p_2(n)$ modulo $5$ and $p_{8}(n)$ modulo $11$. Furthermore,
using the invariance property of $S_{\gamma,\lambda}$ under the Hecke operator $T_{\ell^2}$, we obtain two classes of congruences pertaining to the $m^k$-adic property of $p_r(n)$.



\section{Introduction}
 \label{Sec-Introduction}

 The objective of this paper is to use the theory of modular forms to derive certain congruences of multipartitions modulo powers of primes.

Recall that an ordinary \emph{partition} $\lambda$ of a nonnegative integer $n$ is a non-increasing
sequence of positive integers whose sum is $n$, where $n$ is called the weight of $\lambda$.
The partition function $p(n)$ is defined to be the number of partitions of $n$. A \emph{multipartition} of $n$ with $r$ components, as called by Andrews \cite{Andrews2007},
also referred to as  an $r$-colored partition, see, for example \cite{Brown2009,Eichhorn-Ono1996},  is an
$r$-tuple $\lambda=(\lambda^{(1)},\ldots,\lambda^{(r)})$ of partitions whose weights sum to $n$.
The number of $r$-component multipartitions of $n$ is denoted by $p_r(n)$.

Multipartitions arise
in combinatorics, representation theory and physics. As pointed out by Fayers \cite{Fayers2006}, the
representations of the Ariki-Koike algebra are naturally indexed by multiparititions.
Bouwknegt \cite{Bouwknegt2002} showed that
the Durfee square
formulas of multipartitions are useful in deriving
expressions for the
characters of modules of affine Lie algebras in terms of the universal chiral
partition functions.

For the partition function $p(n)$, Ramanujan \cite{Ramanujan1916,Ramanujan1919,Ramanujan1920,Ramanujan1921}
proved that
\begin{equation}
  \label{eq-Ramanujan-congruence}
 p( An+B )\equiv 0\pmod{M},
\end{equation}
for all nonnegative integers $n$ and for $ (A,B,M)= (5,4,5)$, $ (7,5,7)$  and $(11,6,11)$.
In general, congruences of form  \eqref{eq-Ramanujan-congruence} are called Ramanujan-type congruences. For $m=5$ and $7$,
Watson \cite{Watson1938} proved that
\begin{equation}\label{eq-Ramanujan-5,7,11}
 p( m^k n +\beta_{m,k}) \equiv 0 \pmod{m^k},
\end{equation}
where $k\ge1$ and $\beta_{m,k} \equiv 1/24 \pmod{m^k}$. Atkin \cite{Atkin1967} showed that \eqref{eq-Ramanujan-5,7,11} is also valid for $m=11$. When $M$ is not a power of $5,7$ or $11$,  Atkin and O'Brien \cite{AtkinOBrien1967} discovered the following congruence
\[
 p(11^3\cdot13n+237)\equiv0\pmod{13}.
\]
Using the theory of modular forms, Ono \cite{Ono2000}
proved that for any prime $m\ge 5$ and positive integer $k$, there is a
positive proportion of primes $\ell$ such that
\begin{equation}
 \label{eq-Ono-positiveprop}
 p\left(\frac{m^k \ell^3n+1}{24}\right)\equiv 0 \pmod{m}
\end{equation}
holds for every nonnegative integer $n$  coprime to $\ell$.
 Weaver
\cite{Weaver2001} gave an algorithm for finding the values of $\ell$ in \eqref{eq-Ono-positiveprop} for primes $13\le m\le 31$.

Ramanujan-type congruences of $p_r(n)$ have been extensively studied, see for example \cite{Gandhi1963, Gordon1983,Andrews2007,Atkin1968, Kiming1992, Newman1957, Trenner2006, Gupta1958}. Gandhi \cite{Gandhi1963} derived the following congruences of $p_r(n)$ by applying the identities of Euler and Jacobi
\begin{eqnarray}
 p_2(5n+3)&\equiv& 0 \pmod 5,\label{Gandhi-p2}\\[5pt]
 p_8(11n+4)&\equiv& 0 \pmod {11}.\label{Gandhi-p8}
\end{eqnarray}
With the aid of Sturm's theorem \cite{Sturm1984}, Eichhorn and Ono \cite{Eichhorn-Ono1996} computed an upper bound $C(A,B,r,m^k)$ such that
\[
  p_r(An+B) \equiv 0 \pmod{m^k}
\]
holds for all nonnegative integers $n$ if and only if it is true for $n\le C(A,B,r,m^k)$. For example, to prove \eqref{Gandhi-p2}, it suffices
to check that it holds for $n\leq 3$. In the same vain, one can prove \eqref{Gandhi-p8} by verifying
that it holds for  $n\le 11$.
Treneer \cite{Trenner2006} extended
\eqref{eq-Ono-positiveprop}  to weakly holomorphic modular forms and  showed that for any prime $m \ge 5$ and positive integers $k$, there is a positive proportion of primes $\ell$ such that
\[
p_r\left(\frac{m^k\ell^{\mu_r} n+r}{24}\right) \equiv 0 \pmod{m}
\]
for every nonnegative integer $n$ coprime to $\ell$, where $\mu_r$ equals to
$1$ if $r$ is even and $3$ if $r$ is odd.

The aim of this paper is to study congruence properties of $p_r(n)$ modulo powers of primes. For example, we shall show that
\begin{equation}\label{eq-pr-explicit}
 p_r\left(\frac{m^k\ell^{2\mu K-1}n+r}{24}\right)\equiv 0 \pmod{m^k},
 \end{equation}
 where $r$ is an odd integer, $\ell$ is a prime other than $2,3$ and $m$,  and $\mu$ is a positive integer, $K$ is a fixed positive integer, and $n$ is a positive integer coprime to $\ell$.

To derive congruences of $p_r(n)$, one may consider  the congruence properties of the generating functions of $p_r(n)$. For the case of ordinary partitions, i.e., $r=1$, Chua \cite{Chua2003}
showed that
\begin{equation}\label{eq-Chua-explicit}
\sum_{m n \equiv -1 \rule{-5pt}{0pt} \pmod{24}}  p \left(\frac{mn+1}{24}\right)q^n
 \equiv \eta(24z)^{\gamma_{m}}\phi_{m}(24z) \pmod{m},
\end{equation}
where $\eta(z)$ is Dedekind's eta function, $\gamma_{m}$ is an integer depending on $m$ and $\phi_{m}(z)$ is a
holomorphic modular form.
Ahlgren and Boylan \cite{AhlgrenBoylan2003} extended \eqref{eq-Chua-explicit} to congruences  modulo powers of primes, namely,
 \begin{equation}\label{eq-Ahlgren-explicit}
  F_{m,k}(z) = \sum_{m^k n \equiv -1 \rule{-5pt}{0pt} \pmod{24}}  p \left(\frac{m^kn+1}{24}\right)q^n\equiv \eta(24z)^{\gamma_{m,k}}\phi_{m,k}(24z) \pmod{m^k},
 \end{equation}
 where $\gamma_{m,k}$ is an integer and $\phi_{m,k}(z)$ is a holomorphic modular form.

In order to prove the existence of congruences  of $p_r(n)$ modulo powers of primes,
Brown and Li\cite{Brown2009} introduced the generating function
\begin{equation}\label{eq-BL}
 G_{m,k,r}(z) \equiv \sum_{ \left(\frac{n}{m}\right)= -\left(\frac{-r}{m}\right)}  p_r \left(\frac{n+r}{24}\right)q^n \pmod{m^k},
\end{equation}
and showed that $G_{m,k,r}(z)$ is a modular form of level $576m^3$. Kilbourn \cite{Kilbourn2007} used
the generating function
\begin{equation}\label{eq-Kil}
H_{m,k,r}(z) \equiv \sum_{m n \equiv -r \rule{-5pt}{0pt} \pmod{24}}  p_r \left(\frac{mn+r}{24}\right)q^n \pmod{m^k},
\end{equation}
and proved that $H_{m,k,r}(z)$ is a modular form of level $576m$. However, due to the large
dimensions of the spaces $M_\lambda(\Gamma_0(576m^3))$ and $M_\lambda(\Gamma_0(576m))$,
it does not seem to be a feasible task to compute explicit bases. In other words, to derive explicit congruence formulas of $p_r(n)$, it is
desirable to find a generating function of $p_r(n)$ that can be expressed in terms of  modular forms of a small level.

In this paper, we find the following
extension of the generating function $F_{m,k}(z)$, namely,
\begin{equation}\label{eq-MP-gfun}
F_{m,k,r}(z) = \sum_{m^k n \equiv -r \rule{-5pt}{0pt} \pmod{24}}  p_r \left(\frac{m^kn+r}{24}\right)q^n,
\end{equation}
where $q=e^{2\pi i z}$. We show that
$F_{m,k,r}(z)$ is congruent to
a meromorphic function modulo  $m^k$. More precisely, we find
\begin{equation}
 \label{eq-MP-gfun-intro}
F_{m,k,r}(z) \equiv \eta(24z)^{\gamma_{m,k,r}}\phi_{m,k,r}(24z) \pmod{m^k},
\end{equation}
where $\gamma_{m,k,r}$ is an integer and $\phi_{m,k,r}(z)$ is a holomorphic modular form in $M_{\lambda_{m,k,r}}({\rm SL}_2(\mathbb{Z}))$. Noting that any element of $M_{\lambda_{m,k,r}}({\rm SL}_2 (\mathbb{Z}))$ can be expressed as a polynomial of the Eisenstein series $E_4(z)$ and $E_6(z)$. This enable us to derive explicit congruences of generating functions of $p_r(n)$ modulo $m^k$.

If $\phi_{m,k,r}(z)=0$, then \eqref{eq-MP-gfun-intro} yields a Ramanujan-type congruence as follows
\begin{equation}\label{eq-pr-Ramanujan}
p_r\left(\frac{m^kn+r}{24}\right) \equiv 0 \pmod{m^k}.
\end{equation}
For example, it is easily checked that $\phi_{5,1,2}(z)=0$ and $\phi_{11,1,2}(z)=0$, hence Gandhi's congruences \eqref{Gandhi-p2} and \eqref{Gandhi-p8} are the consequences of \eqref{eq-pr-Ramanujan}.
We also find
\begin{align}
p_2(5^2n+23) &\equiv 0 \pmod{5^2},\\[5pt]
p_8(11^2n+81) &\equiv 0 \pmod{11^2},
\end{align}
since $\phi_{5,2,2}(z)=0$ and $\phi_{11,2,8}(z)=0$.
For more congruences of form \eqref{eq-pr-Ramanujan}, see Table \ref{table-Ramanujan-type}.

On the other hand, if $\phi_{m,k,r}(z)\ne 0$ in \eqref{eq-MP-gfun-intro}, we may use Yang's method \cite{Yang2011} to find congruences of form \eqref{eq-pr-explicit}.
For example, since $F_{5,2,3}(z)$ is congruent to a modular form in the invariant
 space $S_{21,48}$ of $T_{5^2}$ modulo $5^2$, we have
 \[
  p_3\left(\frac{5^2\cdot 13^{199}n+3}{24}\right)\equiv 0\pmod{5^2}.
 \]

\section{Preliminaries}
\label{Sec-Prelim}

To make this paper self-contained, we recall some definitions and facts on modular forms.
In particular, we shall use the $U$-operator, the $V$-operator, the Hecke operator and the twist operator on the modular forms.

Let $k\in\frac{1}{2} \mathbb{Z}$ be an integer or a half-integer, $N$ be a positive integer (with  $4|N$ if $k\not\in\mathbb{Z}$) and $\chi$ be a Nebentypus character. We use $M_k(\Gamma_0(N), \chi)$ to
denote the space of holomorphic modular forms on $\Gamma_0(N)$ of weight $k$ and character $\chi$. The corresponding space of cusp forms is denoted by $S_k(\Gamma_0(N), \chi)$. If $\chi$ is the trivial character, we shall write $M_k(\Gamma_0(N))$ and $S_k(\Gamma_0(N))$ for $M_k(\Gamma_0(N), \chi)$ and $S_k(\Gamma_0(N), \chi)$. Moreover, we write ${\rm SL}_2(\mathbb{Z})$ for $\Gamma_0(1)$.

Let $f(z) \in M_k(\Gamma_0(N), \chi)$ with the following Fourier expansion at $\infty$
\[
f(z)=\sum_{n \ge 0} a(n)q^n,
\]
where $q=e^{2\pi i z}$.
Let us recall some operators acting on $f(z)$.

Let
\[
\gamma = \left( \begin{array}{cc} a & b \\ c & d \end{array} \right)
\]
be a $2 \times 2$ real matrix with positive determinant. The \emph{$k$ slash operator} $|_k$ is defined by
\begin{equation}
 \label{defi-slash-oper}
 (f|_k\gamma)(z) = (\det \gamma)^{k/2} (cz+d)^{-k} f(\gamma z),
\end{equation}
where
\[
\gamma z = \frac{az +b}{cz+d}.
\]
In particular, let $\ell$ be an integer and
\[
\gamma_\ell =  \left(\begin{array}{cc}
     0 & -1 \\
     \ell  & 0
  \end{array}\right).
\]
The \emph{Fricke involution} $W_\ell$ is given by
\begin{equation}\label{defi-W-operator}
f| W_\ell = f|_k \gamma_\ell.
\end{equation}

The \emph{$U$-operator} $U_\ell$ and \emph{$V$-operator} $V_\ell$ are defined by
\begin{equation}\label{defi-U-operator}
f(z)| U_\ell =\sum_{n\ge 0} a(\ell n) q^n,
\end{equation}
and
\begin{equation}\label{defi-V-operator}
f(z)|V_\ell = \sum_{n\ge 0} a(n) q^{\ell n}.
\end{equation}
It is known that
\begin{equation}\label{prop-U-oper}
 f(z) |_k U_\ell = \ell^{\frac{k}{2}-1} \sum_{\mu=0}^{\ell-1} f(z) \Big|_k
 \left(\begin{array}{cc}
  1  & \mu \\
  0  & \ell
 \end{array}\right).
\end{equation}

Let $\psi$ be a Dirichlet character. The $\psi$-twist of $f(z)$
is defined by
\[
 (f \otimes \psi) (z) = \sum_{n\ge 0} \psi(n)a(n)q^n.
\]

Let $\ell$ be a prime and $f(z) \in M_{\lambda+\frac{1}{2}}(\Gamma_0(N), \chi)$ be a modular form of half-integral weight.
The Hecke operator $T_{\ell^2}$ is defined by
\begin{equation}
 \label{defi-Hecke-oper}
f(z)| T_{\ell^2} = \sum_{n \ge 0}
  \left(a(\ell^2 n)+\chi(\ell)\left(\frac{(-1)^\lambda n}{\ell}\right)
  \ell^{\lambda-1}a(n)+\chi(\ell^2) \ell^{2\lambda-1}
  a\left(\frac{n}{\ell^2}\right)\right)q^n.
\end{equation}

We will use the following level reduction properties of the operators $U_\ell$ and $F_\ell=U_\ell + \ell^{\frac{k}{2}-1} W_\ell$ (see \cite[Lemma $1$]{Li1975} and \cite[Lemma $2.2$]{Chua2003}).

\begin{lem}
 \label{lem-oper}
Let $k\in\mathbb{Z}$, $N$ be a positive integer, $\chi$ be a character modulo $N$, and $f(z) \in M_k(\Gamma_0(N), \chi)$.
Assume that  $\ell$ is a prime factor of $N$ and $\chi$ is also a character modulo $N/\ell$.
  \begin{enumerate}
   \item[\rm (1)] If $\ell^2 \mid N$, then $f|U_\ell \in M_k(\Gamma_0(N/\ell), \chi)$.
   \item[\rm (2)] If $N=\ell$ and $\chi$ is the trivial character, then $f|F_\ell \in M_k( {\rm SL}_2(\mathbb{Z}) )$.
  \end{enumerate}
\end{lem}

In the proof of congruence \eqref{eq-MP-gfun-intro} on the generating function $F_{m,k,r}(z)$, we need the following relation
\begin{equation}
 \label{eq-eta-transformation}
  \eta(\gamma z) = \epsilon_{a,b,c,d} (cz+d)^{\frac{1}{2}} \eta(z),
\end{equation}
where
$\gamma=
 \left(\begin{array}{cc}
   a &  b \\
   c &  d
 \end{array}\right)
 \in {\rm SL}_2(\mathbb{Z})$, $\epsilon_{a,b,c,d}$ is a
$24$-th root of unity, and
$\eta(z)$ is Dedekind's eta function as given by
\begin{equation}\label{defi-eta-function}
 \eta(z) = q^{\frac{1}{24}}\prod_{n=1}^\infty (1-q^n).
\end{equation}
As a special case, we have
\begin{equation}
 \label{eq-eta-trans}
 \eta(-1/z)=\sqrt{z/i}\cdot \eta(z).
\end{equation}

\section{Generating functions of $p_r(n)$ modulo $m^k$} \label{Sec-MP-gfun}

In this section, we derive the congruence of the generating function $F_{m,k,r}(z)$ defined by
\eqref{eq-MP-gfun}, namely,
\[
F_{m,k,r}(z) = \sum_{m^k n \equiv -r \rule{-5pt}{0pt} \pmod{24}}  p_r \left(\frac{m^kn+r}{24}\right)q^n.
\]

\begin{thm}\label{thm-MP-gfun}
Let $m\ge 5$ be a prime, and let $k$ and $r$ be positive integers. Then there exists a
modular form $\phi_{m,k,r}(z)\in M_{\lambda_{m,k,r}}({\rm SL}_2 (\mathbb{Z}))$,
such that
\begin{equation}\label{eq-Fphi}
F_{m,k,r}(z) \equiv \eta(24z)^{\gamma_{m,k,r}}\phi_{m,k,r}(24z) \pmod{m^k},
\end{equation}
where
\begin{align}
 \lambda_{m,k,r} & =
  \left\{\begin{array}{ll}
   \frac{m^k-m^{k-1}}{2} r - \frac{\gamma_{m,k,r}+r}{2}, & \mbox{ if $k$ is odd,} \\[7pt]
   (m^k-m^{k-1}) r - \frac{\gamma_{m,k,r}+r}{2}, & \mbox{ if $k$ is even,}
  \end{array}\right. \label{defi-lambda}\\[10pt]
 \gamma_{m,k,r} & = \frac{24\beta_{m,k,r}-r}{m^k},\label{defi-gamma}
\end{align}
and $\beta_{m,k,r}$ is the unique integer in the range $0 \le \beta_{m,k,r}<m^k$ congruent to $r/24$ modulo $m^k$.
\end{thm}

The first step of the proof of Theorem \ref{thm-MP-gfun} is to express $F_{m,k,r}(z)$ in terms of a modular form.
Consider the $\eta$-quotient
\begin{equation}\label{eq-MP-ffun}
 f_{m,k,r}(z) = \left(\frac{\eta(m^kz)^{m^k}}{\eta(z)}\right)^r,
\end{equation}
which is a cusp form in $S_{\frac{(m^k-1)r}{2}}\left(\Gamma_0(m^k),
 \left(\frac{\cdot}{m}\right)^{kr}\right)$.
The following lemma shows that $F_{m,k,r}(z)$ can be obtained from $f_{m,k,r}(z)$ by applying a $U$-operator and a $V$-operator.

\begin{lem}\label{lem-fun-F}
Let $m\ge 5$ be a prime, and let $k$ and $r$ be positive integers. Then we have
\begin{equation}\label{eq-fun-F-operator}
 F_{m,k,r}(z) =\frac{\left(f_{m,k,r}(z)| U_{m^k}\right)\mid V_{24}}
                {\eta(24z)^{m^k r}} .
\end{equation}
\end{lem}

\pf Since
\[
\sum _{n=0}^\infty p_r(n) q^{n}=\prod_{n=1}^\infty\frac{1}{ (1-q^n)^r},
\]
  we find
\begin{eqnarray*}
f_{m,k,r}(z)&=& q^{\frac{m^{2k}-1}{24}r}\prod_{n=1}^\infty\frac{1}{ (1-q^n)^r}\cdot
\prod_{n=1}^\infty(1-q^{m^kn})^{m^k r}\\
&=&  q^{\frac{m^{2k}-1}{24}r}\sum _{n=0}^\infty p_r(n) q^{n}\cdot
\prod_{n=1}^\infty(1-q^{m^kn})^{m^k r}.
\end{eqnarray*}
Applying the operator $U_{m^k}$, we obtain
\[
f_{m,k,r}(z)| U_{m^k}
= \sum _{n=0}^\infty p_r(m^kn+\beta_{m,k,r}) q^{n + \frac{ r (m^{2k}-1) + 24 \beta_{m,k,r}}
{24 m^k} } \cdot \prod_{n=1}^\infty(1-q^{n})^{m^kr},
\]
where $0 \le \beta_{m,k,r}\le m^k-1$ is determined by $24\beta_{m,k,r}\equiv r \pmod{m^k}$.
So we deduce that
\[
\sum _{n=0}^\infty p_r(m^kn+\beta_{m,k,r}) q^{n + \frac{ r (m^{2k}-1) + 24 \beta_{m,k,r}}
{24 m^k}}
= \frac{f_{m,k,r}(z) | U_{m^k} } {\prod_{n=1}^\infty(1-q^{n})^{m^kr}}.
\]
Applying the operator $V_{24}$, we get
\begin{equation}\label{eq-gf-f}
\sum _{n=0}^\infty p_r(m^kn+\beta_{m,k,r}) q^{24n+\frac{24\beta_{m,k,r}-r}{m^k}}
= \frac{\left(f_{m,k,r}(z)| U_{m^k}\right)\mid V_{24}}{\eta(24z)^{m^kr}}.
\end{equation}
Replacing $24n+\frac{24\beta_{m,k,r}-r}{m^k}$
by $n$ in \eqref{eq-gf-f}, or equivalently,
\[
n \to \frac{n}{24} - \frac{24\beta_{m,k,r}-r}{24 m^k},
\]
one sees that the sum on the
left hand side can be written in the form of $F_{m,k,r}(z)$. This completes the
proof. \qed

The second step of the proof of Theorem~\ref{thm-MP-gfun} is to derive a  congruence relation for $f_{m,k,r}(z)|U_{m^k}$ modulo $m^k$.

\begin{thm}
 \label{thm-explict-cong-lem}
Let $m\ge5$ be a prime, and let $k$ and $r$ be positive integers. Then there
exists a modular form
$G_{m,k,r}(z)\in M_{w_{m,k,r}}({\rm SL}_2(\mathbb{Z}))$
such that
\[
f_{m,k,r}(z)|U_{m^k} \equiv G_{m,k,r}(z) \pmod{m^k},
\]
where
\[
 w_{m,k,r} =
 \left\{\begin{array}{ll}
   \frac{2m^k - m^{k-1} -1}{2}r, &\mbox{ if $k$ is odd,}\\[7pt]
   \frac{3m^k - 2m^{k-1} -1}{2}r, &\mbox{ if $k$ is even.}
 \end{array}\right.
 \]
\end{thm}

\noindent {\it Proof.}
Let
\[
g_{m,k,r}(z) = \left(\frac{\eta(z)^m}{\eta(mz)}\right)^{c_k \, m^{k-1}r },
\]
where
\[
c_k = \begin{cases}
    1, & \mbox{if $k$ is odd,} \\[7pt]
    2, & \mbox{if $k$ is even.}
    \end{cases}
\]
 Since $g_{m,k,r}(z)$ is an $\eta$-quotient,
  using the modular transformation property due to  Gordon, Hughes, and Newman \cite{GordonHughes1993,Newman1956,Newman1959}, see also,  \cite[Theorem 1.64]{Ono-CBMS}, we deduce that
\[
g_{m,k,r}(z) \in M_{\frac{c_k\, (m^k-m^{k-1})r}{2}} \left( \Gamma_0(m), \left(\frac{\cdot}{m}\right)^{kr} \right).
\]
Moreover, since $1-q^{mn} \equiv (1-q^n)^m \pmod m$, we see that
\begin{equation}
 \label{eq-g-mod}
 g_{m,k,r}(z) \equiv 1 \pmod{m^k}.
\end{equation}

Since $
f_{m,k,r}(z)  \in S_{\frac{(m^k-1)r}{2}}\left(\Gamma_0(m^k),
\left(\frac{\cdot}{m}\right)^{kr}\right)$,
using Lemma \ref{lem-oper} repeatedly, we obtain that
\[
f_{m,k,r}(z)|U_{m^{k-1}}
\in S_{\frac{(m^k-1)r}{2}} \left( \Gamma_0(m), \left(\frac{\cdot}{m}\right)^{kr} \right).
\]
Thus, $f_{m,k,r}(z)|U_{m^{k-1}} \cdot g_{m,k,r}(z)$
is a modular form on $\Gamma_0(m)$ of the trivial character and of weight
\[
w_{m,k,r} = \frac{c_k\, (m^k-m^{k-1})r}{2} + \frac{(m^k-1)r}{2}.
\]
Invoking Lemma \ref{lem-oper}, we find that
 \begin{equation}
  \label{eq-fg-Foper}
  G_{m,k,r}(z) = \left(f_{m,k,r}(z)|U_{m^{k-1}} \cdot g_{m,k,r}(z)\right)| F_m
 \end{equation}
is a modular form in
$M_{w_{m,k,r}}({\rm SL}_2(\mathbb{Z}))$. 

To complete the proof of Theorem~\ref{thm-explict-cong-lem}, it remains to show that 
\begin{equation}\label{eq-f-F}
\left(f_{m,k,r}(z)|U_{m^{k-1}} \cdot g_{m,k,r}(z)\right)| F_m \equiv f_{m,k,r}(z)|U_{m^{k}}  \pmod{m^k},
\end{equation}
where
\[
F_m=U_m + m^{\frac{w_{m,k,r}}{2}-1} W_m,
\]
and the operator $W_m$ is given by \eqref{defi-W-operator}.
By congruence \eqref{eq-g-mod}, we see that the left hand side of \eqref{eq-f-F} equals
\[
  f_{m,k,r}(z)|U_{m^k}  + m^{\frac{w_{m,k,r}}{2}-1} \left(f_{m,k,r}(z)|U_{m^{k-1}} \cdot g_{m,k,r}(z)\right)|W_m \pmod{m^k}.
 \]

To prove \eqref{eq-f-F}, it suffices to show that
\begin{equation}\label{eq-fg-Woper}
 m^{\frac{w_{m,k,r}}{2}-1} \left(f_{m,k,r}(z)|U_{m^{k-1}}
    \cdot g_{m,k,r}(z)\right)|W_m \equiv 0\pmod{m^k}.
\end{equation}
We only consider the case when $k$ is odd. The case when $k$ is even can be dealt with in the same manner. In light of the transformation formula \eqref{eq-eta-trans} of the
eta function, we find that
\begin{eqnarray*}
g_{m,k,r}(z)| W_m &=&  m^{\frac{(m^k-m^{k-1})r}{4}}
   (mz)^{-\frac{(m^k-m^{k-1})r}{2}}g_{m,k,r}
   \left(-\frac{1}{mz}\right) \\[7pt]
 &=& m^{-\frac{(m^k-m^{k-1})r}{4}}z^{-\frac{(m^k-m^{k-1})r}{2}}
 \left(\frac{(\sqrt{mz/i}\, \eta(mz))^m}{\sqrt{z/i}\, \eta(z)}\right)^{m^{k-1}r}\\[7pt]
 &=&m^{\frac{(m+1)m^{k-1}r}{4}}(-i)^{\frac{(m-1)m^{k-1}r}{2}}
 \left(\frac{\eta(mz)^m}{\eta(z)}\right)^{m^{k-1}r}.
\end{eqnarray*}
Therefore,   \eqref{eq-fg-Woper}  can be deduced from the following congruence
\begin{equation}\label{eq-UW}
  m^{\frac{(3m^k-1)r}{4}-1}\left.(\left. f_{m,k,r}(z)
   \right| U_{m^{k-1}}) \right| W_m \equiv 0 \pmod{m^k}.
\end{equation}
By the property of $U$-operator as in \eqref{prop-U-oper}, we have
\begin{eqnarray}
&& m^{\frac{(3m^k-1)r}{4}-1}\left.\left.f_{m,k,r}(z)\right|
  U_{m^{k-1}}\right|W_m\nonumber \\[7pt]
&=&m^{\frac{(k+2)m^kr-(r+4)k}{4}}\sum_{\mu=0}^{m^{k-1}-1}
   f_{m,k,r}(z)
   \Big|_{\frac{(m^k-1)r}{2}}
\left(\begin{array}{cc}
 1   &   \mu    \\
 0   &   m^{k-1}
\end{array}\right)\Big|W_m\nonumber\\[7pt]
&=&m^{\frac{(k+2)m^kr-(r+4)k}{4}}
   \sum_{\mu=0}^{m^{k-1}-1}
   \left.f_{m,k,r}(z)
   \right|_{\frac{(m^k-1)r}{2}}
\left(\begin{array}{cc}
 \mu m   & -1  \\
 m^k   &  0
\end{array}\right). \label{eq-involution-modulo0}
\end{eqnarray}
Using the transformation formula \eqref{eq-eta-trans} of the eta function, \eqref{eq-involution-modulo0} can be written as
\begin{eqnarray}
  &&m^{\frac{m^kr}{2} - k} z^{-\frac{(m^k-1)r}{2}}
  \sum_{\mu=0}^{m^{k-1}-1}
   \left(\frac{\eta(m\mu-\frac{1}{z})^{m^k}}{\eta(\frac{m\mu z-1}{m^kz})}\right)^r
 \nonumber\\[7pt]
&=&m^{\frac{m^kr}{2} - k} z^{\frac{r}{2}}
  \eta(z)^{m^kr}\sum_{\mu=0}^{m^{k-1}-1}
\frac{\alpha_\mu}{\eta(\frac{ m\mu z-1}{m^kz})^r},
\label{eq-involution-modulo}
\end{eqnarray}
where $\alpha_\mu$ is a certain $24$-th root unity.

For $\mu \neq 0$, we write
$\mu = m^s t$ where $m \!\nmid\! t$.
For $\mu=0$, we set $s=k-1$ and $t=0$. In either case,  there exist integers $b$ and $d$ such that
$bt+dm^{k-s-1}=-1$. It follows that
\[
\left(\begin{array}{cc}
m\mu  &   -1    \\
 m^k  &   0
\end{array}\right)
=
\left(\begin{array}{cc}
 t  &   d    \\
 m^{k-s-1} &   -b
\end{array}\right)
\left(\begin{array}{cc}
 m^{s+1}  &   b    \\
 0  &   m^{k-s-1}
\end{array}\right).
\]
Applying the corresponding slash operator to $\eta(z)$, we obtain that
\[
 \eta\left(\frac{m\mu z-1}{m^kz}\right)=
 \epsilon_\mu m^{\frac{s+1}{2}}z^{\frac{1}{2}}
 \eta\left(\frac{m^{s+1}z+b}{m^{k-s-1}}\right),
\]
where $\epsilon_\mu$ is a $24$-th root of unity. Since the coefficients of the Fourier expansion of $\eta(z)$ at $\infty$ are integers and the coefficient of the term with the lowest degree is $1$, the Fourier coefficients of each term in \eqref{eq-involution-modulo} are divisible by $m^{\frac{m^k-s-1}{2}r - k}$ in the ring $\mathbb{Z}[\zeta_{24}]$.
Clearly,  $0 \le s \le k-1$. Thus we have
\[
\frac{m^k-s-1}{2}r - k \ge \frac{m^k-k}{2}r - k \ge \frac{m^k-k}{2} - k \ge k,
\]
for $m \ge 5$ and $k \ge 1$. Hence the Fourier coefficients of each term in \eqref{eq-involution-modulo} are divisible by $m^k$. So we arrive at \eqref{eq-UW}. This completes the proof. \qed

 We are now in a position to finish the proof of Theorem~\ref{thm-MP-gfun}.

\noindent {\it Proof of Theorem~\ref{thm-MP-gfun}.}
By Theorem \ref{thm-explict-cong-lem}, there exists a modular form
$G_{m,k,r}(z)\in M_{w_{m,k,r}}({\rm SL}_2(\mathbb{Z}))$
such that
\begin{equation}
 \label{eq-f-G}
 f_{m,k,r}(z)|U_{m^k} \equiv G_{m,k,r}(z) \pmod{m^k}.
\end{equation}
Let
\[
\phi_{m,k,r}(z) = \frac{G_{m,k,r}(z)}{\Delta(z)^{\frac{m^k r + \gamma_{m,k,r}}{24}}},
\]
where $\Delta(z)=\eta(z)^{24}$ is Ramanujan's $\Delta$-function.
In the proof of Lemma~\ref{lem-fun-F} we have shown that
\[ f_{m,k,r}(z)| U_{m^k}
= \sum _{n=0}^\infty p_r(m^kn+\beta_{m,k,r}) q^{n + \frac{ r (m^{2k}-1) + 24 \beta_{m,k,r}}
{24 m^k} } \cdot \prod_{n=1}^\infty(1-q^{n})^{m^kr},
\]
which implies that the order of the Fourier expansion of $f_{m,k,r}(z)| U_{m^k}$ at $\infty$ is at least
\[
\frac{r (m^{2k}-1) + 24 \beta_{m,k,r}} {24 m^k} = \frac{m^k r + \gamma_{m,k,r}}{24}.
\]
Thus $\phi_{m,k,r}(z)$ is a modular form in $M_{\lambda_{m,k,r}}({\rm SL}_2(\mathbb{Z}))$.
Combining \eqref{eq-f-G} and Lemma \ref{lem-fun-F}, we conclude that
\begin{align*}
 F_{m,k,r}(z)
 & \equiv  \frac{\left.\left(\Delta(z)^{\frac{m^k r + \gamma_{m,k,r}}{24}} \phi_{m,k,r}(z)\right)\right|V_{24}}{\eta(24z)^{m^kr}}\\[7pt]
 & = \eta(24z)^{\gamma_{m,k,r}}\phi_{m,k,r}(24z) \pmod{m^k},
\end{align*}
as required. \qed

\section{Congruences of $p_r(n)$ modulo $m^k$}
\label{Sec-MP-cong}

In this section, we apply Theorem~\ref{thm-MP-gfun} on the congruence relation for the generating function $F_{m,r,k}(z)$ and Yang's method \cite{Yang2011} to derive two classes of congruences of $p_r(n)$ modulo $m^k$.

Let
\[
S_{\gamma, \lambda} = \{ \eta(24 z)^\gamma \phi(24 z) \colon \phi(z) \in M_\lambda({\rm SL}_2(\mathbb{Z})) \}.
\]
Yang \cite{Yang2011} showed that when $\gamma$ is an odd integer such that $0<\gamma<24$ and $\lambda$ is a nonnegative even integer, $S_{\gamma,\lambda}$ is an invariant subspace of $S_{\lambda+\gamma/2}(\Gamma_0(576), \chi_{12})$ under the action of the Hecke algebra. More precisely, for all primes $\ell \not= 2,3$ and all $f \in S_{\gamma,\lambda}$, we have $f|T_{\ell^2} \in S_{\gamma,\lambda}$.
By the invariant property of $S_{\gamma,\lambda}$,    we obtain two classes of congruences of $p_r(n)$ modulo $m^k$.

\begin{thm}
 \label{thm-MP-cong}
 Let $m\ge5$ be a prime, $k$ be a positive integer, $r$ be an odd positive integer less than $m^k$, and $\ell$ be a prime different from
  $2,3$ and $m$. Then
 there exists an explicitly computable positive integer $K$ such that
 \begin{equation}
  \label{eq-MP-explicitcong}
 p_r\left(\frac{m^k\ell^{2\mu K-1}n+r}{24}\right)\equiv 0 \pmod{m^k}
 \end{equation}
for all positive integers $\mu$ and all positive integers $n$ relatively prime to $\ell$.
There is also a positive integer $M$ such that
 \begin{equation}
  \label{eq-MP-ladic-cong}
 p_r\left(\frac{m^k\ell^in+r}{24}\right)\equiv
   p_r\left(\frac{m^k\ell^{2M+i}n+r}{24}\right) \pmod{m^k}
 \end{equation}
 for all nonnegative integers $i$ and $n$.
\end{thm}

\pf
According to congruence relation \eqref{eq-Fphi}, the generating function $F_{m,k,r}(z)$  is congruent to a modular form in $S_{\gamma_{m,k,r},\lambda_{m,k,r}}$, where $\lambda_{m,k,r}$ and $\gamma_{m,k,r}$ are integers as given in \eqref{defi-lambda} and \eqref{defi-gamma}. Let $\{f_1(z),\ldots,f_d(z)\}$
be a $\mathbb{Z}$-basis of the space
$ S_{\gamma_{m,k,r},\lambda_{m,k,r}}\cap \mathbb{Z}[[q]]$ and
\[
 f_i(z) = \sum_{n\ge0} a_i(n) q^n,
\]
where $i=1,\ldots, d$ and $q= e^{2\pi iz}$.

To prove \eqref{eq-MP-explicitcong}, it suffices  to show that there exists a positive integer $K$ such that
\begin{equation}\label{eq-ai}
 a_i\left(\frac{m^k\ell^{2\mu K-1}n+r}{24}\right)\equiv 0 \pmod{m^k}
\end{equation}
for all $n$ coprime to $\ell$ and $i=1,\ldots, d$.

From the relation $ \gamma_{m,k,r}m^k   = 24\beta_{m,k,r}-r$, one sees that
$\gamma_{m,k,r}$ and $r$ have the same parity.
Since $r<m^k$ is odd, we have $0<\gamma_{m,k,r}<24$, and hence $S_{\gamma_{m,k,r},\lambda_{m,k,r}}$ is  invariant under the Hecke operator $T_{\ell^2}$. So there exists a $d\times d$ matrix $A$ such that
\begin{equation}\label{eq-matrix-A}
 \left.
 \left(\begin{array}{c}
  f_1 \\  \vdots \\ f_d
 \end{array}\right)\right|T_{\ell^2}
 = A
  \left(\begin{array}{c}
  f_1 \\  \vdots \\ f_d
 \end{array}\right).
\end{equation}
Let
\[
 X=\left(\begin{array}{cc}
  A &I_d \\
  -\ell^{\gamma_{m,k,r}+2\lambda_{m,k,r}-2}I_d & 0
 \end{array}\right).
\]
Using the property of the basis $\{f_1(z),\ldots,f_d(z)\}$ under the action of the $U$-operator as given by Yang \cite[Corollary 3.4]{Yang2011}, we obtain
\begin{equation}\label{eq-basis-Uoper}
 \left.
 \left(\begin{array}{c}
  f_1 \\  \vdots \\ f_d
 \end{array}\right)\right|U_{\ell^2}^s
 = A_s
  \left(\begin{array}{c}
  f_1 \\  \vdots \\ f_d
 \end{array}\right)
 + B_s
  \left(\begin{array}{c}
  g_1 \\  \vdots \\ g_d
 \end{array}\right)
 +C_s
  \left.\left(\begin{array}{c}
  f_1 \\  \vdots \\ f_d
 \end{array}\right)
 \right|V_{\ell^2} ,
\end{equation}
where $s$ is a positive integer, $g_i=f_i\otimes \left(\frac{\cdot}{\ell}\right)$, and
 $A_s, B_s$ and $C_s$ are $d\times d$ matrices given by
\begin{align}
 &\left(\begin{array}{cc}
  A_s & A_{s-1}
 \end{array}\right)
 =
 \left(\begin{array}{cc}
  I_d & 0
 \end{array}\right)
  X^s,\label{defi-A-matrix}\\[7pt]
 & B_s =  -\ell^{\lambda_{m,k,r}+(\gamma_{m,k,r}-3)/2}
  \left(\frac{(-1)^{(\gamma_{m,k,r}-1)/2}12}{\ell}\right)
  A_{s-1},\nonumber\\[7pt]
 &C_s = -\ell^{\gamma_{m,k,r}+2\lambda_{m,k,r}-2} A_{s-1}.\nonumber
\end{align}
Since $\gcd(m,\ell)=1$, the matrix $X \pmod{m^k}$ is invertible in the ring ${\cal M}$ consisting of $2d\times 2d$ matrices over $\mathbb{Z}_{m^k}$. By the finiteness of ${\cal M}$, we see that there exist integers $a>b$ such that $X^a$ and $X^b$ are linear dependent over $\mathbb{Z}_{m^k}$, i.e., there exists a constant $c \in \mathbb{Z}_{m^k}$ such that $X^a \equiv cX^b \pmod{m^k}$. Thus
$X^K \equiv c I_{2d} \pmod{m^k}$, where $K=a-b$. In view of the relation
\[
 \left(\begin{array}{cc}
  A_{\mu K-1} & A_{\mu K-2}
 \end{array}\right)
 \equiv c^\mu
 \left(\begin{array}{cc}
  I_d & 0
 \end{array}\right)
 X^{-1} \pmod{m^k},
\]
we find that $A_{\mu K-1} \equiv 0 \pmod{m^k}$.
Hence, from \eqref{eq-basis-Uoper} it follows that
\[ \left.
 \left(\begin{array}{c}
  f_1 \\  \vdots \\ f_d
 \end{array}\right)\right|U_{\ell^2}^{\mu K-1}
 \equiv B_{\mu K-1}
  \left(\begin{array}{c}
  g_1 \\  \vdots \\ g_d
 \end{array}\right)
 +C_{\mu K-1}
  \left.\left(\begin{array}{c}
  f_1 \\  \vdots \\ f_d
 \end{array}\right)
 \right|V_{\ell^2} \pmod{m^k}.
\]
Applying the $U$-operator $U_\ell$ to both sides and observing that
\[
g_i | U_\ell = f_i \otimes \left. \left( \frac{\cdot}{\ell} \right) \right| U_\ell = 0,
\]
the relation  \eqref{eq-basis-Uoper}  leads to the
following congruence
\[
 \left.
 \left(\begin{array}{c}
  f_1 \\  \vdots \\ f_d
 \end{array}\right)\right|U_{\ell^2}^{\mu K-1}U_\ell
 \equiv C_{\mu K-1}
  \left.\left(\begin{array}{c}
  f_1 \\  \vdots \\ f_d
 \end{array}\right)
 \right|V_{\ell} \pmod{m^k},
\]
which implies \eqref{eq-ai}.

We now turn to the proof of congruence \eqref{eq-MP-ladic-cong}. By the finiteness of ${\cal M}$, we see that there exists a positive integer $M$ such that
$X^M \equiv I_{2d} \pmod{m^k}$. Thus matrix equation \eqref{defi-A-matrix} reduces to
the following congruence
\[
 \left(\begin{array}{cc}
  A_M & A_{M-1}
 \end{array}\right)
 \equiv
 \left(\begin{array}{cc}
  I_d & 0
 \end{array}\right)  \pmod{m^k}.
\]
It follows that
$ A_M\equiv I_d \pmod{m^k}$ and
$ B_M\equiv C_M\equiv 0\pmod{m^k}$.
Thus, relation \eqref{eq-basis-Uoper} implies
\[
 \left.
 \left(\begin{array}{c}
  f_1 \\  \vdots \\ f_d
 \end{array}\right)\right|U_{\ell^2}^{M}
  \equiv
  \left(\begin{array}{c}
  f_1 \\  \vdots \\ f_d
 \end{array}\right) \pmod{m^k}.
\]
So the  coefficient  of $q^n$ 
is congruent to the coefficient of
$q^{\ell^{2M}n}$¡¡in $f_i(z)$ modulo $m^k$ for all $i$ and $n$. Since $F_{m,k,r}(z)$ is a linear combination of $f_i(z)$ with integer coefficients, we obtain congruence \eqref{eq-MP-ladic-cong}. This completes the proof. \qed

\section{Examples}\label{Sec-examples}

In this section, we present some consequences of Theorem \ref{thm-MP-gfun} and
Theorem \ref{thm-MP-cong}.
We first give some examples for the congruences of the generating function $F_{m,k,r}(z)$ of $p_r(n)$.

\begin{exam} By Theorem \ref{thm-MP-gfun}, we find
\[
 F_{m,k,r}(z) \equiv \eta(24z)^{\gamma_{m,k,r}}\phi_{m,k,r}(24z) \pmod{m^k},
\]
where $\gamma_{m,k,r}$ is an integer,  $\phi_{m,k,r}(z)$ is a polynomial  of $\Delta(z)$ and the Eisenstein series $E_4(z)$ and $E_6(z)$.  Table \ref{table-explicit-cong} gives the list of explicit expressions of  $\eta(z)^{\gamma_{m,1,r}}\phi_{m,1,r}(z)$ for $m \le 19$ and $2 \le r \le 7$.
\end{exam}

\begin{longtable}{lr|l}
  \hline
   $r$ & $m$ & $\eta(z)^{\gamma_{m,1,r}}\phi_{m,1,r}(z)$ \\
   \hline
   $2$ & $5$ & $0$ \\
       & $7$ & $3\eta(z)^{10}$ \\
       & $11$ & $2\eta(z)^2E_4(z)^2$\\
       & $13$ & $8\eta(z)^{22}$ \\
       & $17$ & $5\eta(z)^{14}E_4(z)^2$\\
       & $19$ & $\eta(z)^{10}(14E_4(z)^3+12\Delta(z))$ \\
  \hline
   $3$ & $5$ & $4\eta(z)^9$ \\
       & $7$ & $3\eta(z)^3E_6(z)$ \\
       & $11$ & $0$\\
       & $13$ & $\eta(z)^9(4E_4(z)^3+6\Delta(z))$ \\
       & $17$ & $0$\\
       & $19$ & $\eta(z)^{15}(2E_6(z)^3 +3E_6(z)\Delta(z))$ \\
  \hline
   $4$ & $5$ & $4\eta(z)^4E_4(z)$ \\
       & $7$ & $0$ \\
       & $11$ & $\eta(z)^4(3E_4(z)^4 +8E_4(z)\Delta(z))$\\
       & $13$ & $\eta(z)^{20}(7E_4(z)^3+4\Delta(z))$ \\
       & $17$ & $\eta(z)^4(6 E_4(z)^7 +11E_4(z)^4\Delta(z)+4E_4(z)\Delta(z)^2)$\\
       & $19$ & $\eta(z)^{20}(16E_4(z)^6+18E_4(z)^3\Delta(z)+2\Delta(z)^2)$ \\
  \hline
   $5$ & $5$ & $\eta(z)^{-1}E_4(z)^2$ \\
       & $7$ & $\eta(z)^{13}E_6(z)$ \\
       & $11$ & $0$ \\
       & $13$ & $\eta(z)^7(8E_4(z)^6+11E_4(z)^3\Delta(z)+5\Delta(z)^2)$ \\
       & $17$ & $\eta(z)^{11}(16E_4(z)^8 +16E_4(z)^5\Delta(z)+4E_4(z)^2\Delta(z)^2)$ \\
       & $19$ & $\eta(z)(5E_6(z)^7+15E_6(z)^5\Delta(z)+16E_6(z)^3\Delta(z)^2)$ \\
  \hline
   $6$ & $5$ & $0$ \\
       & $7$ & $\eta(z)^6(6 E_4(z)^3+6\Delta(z))$ \\
       & $11$ & $\eta(z)^6(10E_4(z)^6+E_4(z)^3\Delta(z))$\\
       & $13$ & $\eta(z)^{18}(7E_4(z)^6+8E_4(z)^3\Delta(z)+6\Delta(z)^2)$ \\
       & $17$ & $\eta(z)^{18}(3E_4(z)^9+3E_4(z)^6\Delta(z)+5E_4(z)^3\Delta(z)^2)$\\
       & $19$ & $\eta(z)^6(6E_4(z)^{12}+E_4(z)^{9}\Delta(z)+14\Delta(z)^4)$ \\
  \hline
   $7$ & $5$ & $0$ \\
       & $7$ & $\eta(z)^{-1}E_6(z)^3$ \\
       & $11$ & $0$ \\
       & $13$ & $\eta(z)^5(10E_4(z)^9+6E_4(z)^6\Delta(z)+9E_4(z)^3\Delta(z)^2+11\Delta(z)^3)$ \\
       & $17$ & $\eta(z)(7E_4(z)^{13}+2E_4(z)^{10}\Delta(z)+E_4(z)^7\Delta(z)^2+3E_4(z)^4\Delta(z)^3)$\\
       & $19$ & $0$\\
  \hline
  \\
\caption{\small Explicit congruences derived from Theorem \ref{thm-MP-gfun}.\label{table-explicit-cong}}
\end{longtable}

\begin{exam}
Let $0 \le \beta <m^k$ be an integer with $\beta \equiv r/24 \pmod{m^k}$.
If $\phi_{m,k,r}(z) \equiv 0 \pmod{m^k}$,
 using Theorem \ref{thm-MP-gfun} we obtain the following Ramanujan-type congruences of  multipartition functions
\begin{equation}\label{eq-Rtype}
 p_r( m^kn +\beta ) \equiv 0 \pmod{m^k}.
\end{equation}
The values of $m$ and $\beta$ for
$r\le 9$ and $k=1,2$  are given  in Table \ref{table-Ramanujan-type}.
\end{exam}

\begin{table}[htb]
 \centering
\begin{tabular}{|c|l|l|}
  \hline
  $r$ & $(m,\beta)$ & $(m^2,\beta)$ \\
  \hline
  $1$ & $(5,4),(7,5),(11,6)$ & $(25,24),(49,47),(121,116)$ \\
  \hline
  $2$ & $(5,3)$ & $(25,23)$ \\
  \hline
  $3$ & $(11,7),(17,15)$ & $(121,106)$ \\
  \hline
  $4$ & $(7,6)$ & $(49,41)$ \\
  \hline
  $5$ & $(11,8),(23,5)$ & $(121,96)$ \\
  \hline
  $6$ & $(5,4)$ & $(25,19)$ \\
  \hline
  $7$ & $(5,3),(11,9),(19,9)$ & $(25,18), (121,86)$ \\
  \hline
  $8$ & $(7,5),(11,4)$ & $(121,81)$ \\
  \hline
  $9$ & $(17,11),(19,17),(23,9)$ & \ \ $\rule{18pt}{0.5pt}$ \\
  \hline
\end{tabular}
 \caption{ \small Ramanujan-type congruences of multipartitions.}
 \label{table-Ramanujan-type}
\end{table}

It can be seen that Table~\ref{table-Ramanujan-type} contains  the Ramanujan
congruences \eqref{eq-Ramanujan-congruence} of $p(n)$ modulo $5, 7$ and $11$, as well as Gandhi's congruences \eqref{Gandhi-p2} for $p_2(n)$ and \eqref{Gandhi-p8} for $p_8(n)$.

The following examples demonstrate
how to derive  certain congruences  of $p_r(n)$ with the aid of Theorem \ref{thm-MP-cong}.

\begin{exam} For the values of $\ell$ and $K_\ell$ as given in Table~\ref{table-Yang-3congruence}, we have
 \begin{equation}\label{eq-Yangp3}
  p_3\left(\frac{7\cdot \ell^{2\mu K_\ell-1}n+3}{24}\right)\equiv 0\pmod{7}
 \end{equation}
for all positive integers $\mu$ and all positive integers $n$ not divisible by $\ell$.

\begin{table}[ht]
\centering
\begin{tabular}{r|rrrrrrrrrrrrrr}
\hline
 $\ell$ & $5$ & $11$ & $13$ & $17$ & $19$ & $23$ & $29$ & $31$ & $37$ & $41$ & $43$ & $47$ & $53$ & $59$ \\[5pt]
 \hline
 $a_\ell$ & $6$ & $4$ & $0$ & $4$ & $3$ & $6$ & $2$ & $5$ & $3$ & $0$ & $0$ & $3$ & $5$ & $5$  \\[5pt]
 $K_\ell$ & $6$ & $7$ & $2$ & $6$ & $8$ & $7$ & $7$ & $8$ & $3$ & $2$ & $2$ & $8$ & $3$ & $8$  \\
 \hline
\end{tabular}
\caption{Eigenvalues $a_\ell$ of $F_{7,1,3}(z)$ acted by $T_{\ell^2}$ and the corresponding $K_\ell$.}
\label{table-Yang-3congruence}
\end{table}
\end{exam}

\pf By Theorem
\ref{thm-MP-gfun} we find
\[
 F_{7,1,3}(z) \equiv 3\eta(24z)^3E_6(24z) \pmod{7}.
\]
Since $\eta(24z)^3E_6(24z)$ belongs to the $1$-dimensional space $S_{3,6}$, for any prime $\ell \not=2,3,7$,  there exists an integer $a_\ell$ such that
\[
F_{7,1,3}(z) | T_{\ell^2} \equiv a_\ell F_{7,1,3}(z) \pmod{7}.
\]
 Inspecting the proof of Theorem~\ref{thm-MP-cong}, we obtain the corresponding orders $K_\ell$ for which
 congruence \eqref{eq-Yangp3} holds. \qed

\begin{exam}
We have
 \[
  p_3\left(\frac{5^2\cdot 13^{199}n+3}{24}\right)\equiv 0\pmod{5^2}
 \]
for all integers $n$ coprime to $13$ and
 \[
   p_3\left(\frac{5^2\cdot 13^in+3}{24}\right)\equiv
   p_3\left(\frac{5^2\cdot 13^{200+i}n+3}{24}\right)\pmod{5^2}
 \]
for all nonnegative integers $n$ and $i$.
\end{exam}

\pf
 By Theorem \ref{thm-MP-gfun}, $F_{5,2,3}(z)$ is congruent to a modular form in the
 space $S_{21,48}$ of dimension $5$.
Setting
 \[
  f_i = \eta(24z)^{21} E_4(24z)^{3(5-i)} \Delta(24z)^{i-1},
 \]
for $1\leq i \leq 5$.
Clearly, $f_1, f_2, \ldots , f_5$  form
 a  $\mathbb{Z}$-basis of $S_{21,48}\cap\mathbb{Z}[[q]]$. Let $A$ be the matrix of $T_{\ell^2}$ with
 respect to this basis. By computing the first five Fourier coefficients of $f_i$ and $f_i| T_{13^2}$ and equating the Fourier coefficients of both sides of \eqref{eq-matrix-A}, we find
  \[
  A \equiv \left(\begin{array}{ccccc}
       17& 21& 18& 3& 3\\
       0& 19& 5& 5& 5\\
       0& 0& 22& 4& 19\\
       0& 0& 0& 22& 10\\
       0& 0& 0& 0& 12
      \end{array}\right) \quad \pmod{5^2},
 \]
with the corresponding orders $K=M=100$.  Setting $\mu=1$ in Theorem~\ref{thm-MP-cong},
we complete the proof. \qed

Below are two more examples for $p_3(n)$ and $p_5(n)$ modulo $7^2$. The proofs are analogous  to the proof of the above example, and hence are omitted.

\begin{exam}
We have
 \[
  p_3\left(\frac{7^2\cdot 11^{2351}n+3}{24}\right)\equiv 0\pmod{7^2}
 \]
for all positive integers $n$ coprime to $7$ and
 \[
  p_3\left(\frac{7^2\cdot11^in+3}{24}\right)\equiv
   p_3\left(\frac{7^2\cdot11^{1176+i}n+3}{24}\right) \pmod{7^2}
 \]
for all nonnegative integers $n$ and $i$.
\end{exam}

\begin{exam}
We have
 \[
  p_5\left(\frac{7^2\cdot 17^{195}n+5}{24}\right)\equiv 0\pmod{7^2}
 \]
  for all positive integers $n$ coprime to $17$ and
 \[
  p_5\left(\frac{7^2\cdot17^in+5}{24}\right)\equiv
   p_5\left(\frac{7^2\cdot17^{588+i}n+5}{24}\right) \pmod{7^2}
 \]
 for all nonnegative integers $n$ and $i$.
\end{exam}

\vspace{.2cm} \noindent{\bf Acknowledgments.}
 This work was
supported by the 973 Project, the PCSIRT Project of the Ministry of
Education,  and the National
Science Foundation of China.

\end{document}